\documentclass[11pt]{article}
\usepackage{geometry}                
\usepackage{pgf,tikz}
\usetikzlibrary{arrows}
\usepackage{graphicx}
\usepackage{amssymb,amsbsy,amsmath,amscd,amsthm}
\usepackage{epstopdf}
\usepackage{tikz}
\usetikzlibrary{arrows}

\usepackage[applemac]{inputenc}

\usepackage[english,francais]{babel}

\newtheorem{definition}{Definition}
\newtheorem{theorem}{Theorem}

\newtheorem{lemma}{Lemma}
\newtheorem{remark}{Remark}

\newtheorem{corollary}{Corollary}
\newtheorem*{AssumptionA}{Assumption A}

\newcommand{\R}{\mathbb R}
\newcommand{\PP}{\mathbb P}

\definecolor{ffqqqq}{rgb}{1,0,0}
\definecolor{qqwuqq}{rgb}{0,0.39,0}
\definecolor{xdxdff}{rgb}{0.49,0.49,1}
\definecolor{uququq}{rgb}{0.25,0.25,0.25}
\definecolor{zzttqq}{rgb}{0.6,0.2,0}
\definecolor{cqcqcq}{rgb}{0.75,0.75,0.75}

\title{Adaptive Estimation of Convex Sets and Convex Polytopes \\ from Noisy Data}

\author{Victor-Emmanuel Brunel}

\date{}                                           

\begin{document}

\maketitle{}

\selectlanguage{english}
\begin{abstract} We estimate convex polytopes and general convex sets in $\mathbb R^d,d\geq 2$ in the regression framework. We measure the risk of our estimators using a $L^1$-type loss function and prove upper bounds on these risks. We show that, in the case of polytopes, these estimators achieve the minimax rate. For polytopes, this minimax rate is $\frac{\ln n}{n}$, which differs from the parametric rate for non-regular families by a logarithmic factor, and we show that this extra factor is essential. Using polytopal approximations we extend our results to general convex sets, and we achieve the minimax rate up to a logarithmic factor. In addition we provide an estimator that is adaptive with respect to the number of vertices of the unknown polytope, and we prove that this estimator is optimal in all classes of polytopes with a given number of vertices.

 \end{abstract}
\selectlanguage{french}

\section{Introduction}
\subsection{Definitions and notations}
Let $d\geq 2$ be a positive integer. Assume that we observe a sample of $n$ i.i.d. pairs $(X_i,Y_i), i=1,\ldots, n$ such that $X_1,\ldots,X_n$ have the uniform distribution on $[0,1]^d$ and 
\begin{equation}
\label{posprob} Y_i=I(X_i\in G)+\xi_i, i=1,\ldots,n.
\end{equation} 
The collection $X_1,\ldots,X_n$ is called the design. The error terms $\xi_i, i=1,\ldots,n,$ are i.i.d. random variables independent of the design, $G$ is a subset of $[0,1]^d$, and $I(\cdot\in G)$ stands for the indicator function of the set $G$. Here we aim to estimate the set $G$ in Model \eqref{posprob}.

A subset $\hat G_n$ of $[0,1]^d$ is called a set estimator, or simply, in our framework, an estimator, if it is a Borel set and if there exists a real measurable function $f$ defined on $([0,1]^d\times\mathbb R)^n$ such that $I(\cdot\in \hat G_n)=f(\cdot, X_1,Y_1, \ldots, X_n,Y_n)$.

If $G$ is a measurable (with respect to the Lebesgue measure on $\R^d$) subset of $[0,1]^d$, we denote by $|G|_d$ or, when there is no possible confusion, simply by $|G|$, its Lebesgue measure and by $\mathbb P_G$ the probability measure with respect to the distribution of the collection of $n$ pairs $(X_i,Y_i), i=1,\ldots, n$. Where it is necessary to indicate the dependence on $n$ we use the notation $\mathbb P_G^{\otimes n}$. If $G_1$ and $G_2$ are two measurable subsets of $\mathbb R^d$ their Nikodym pseudo distance $d_1(G_1,G_2)$ is defined as 
\begin{equation}
	\label{defdist} d_1(G_1,G_2)=|G_1\triangle G_2|.
\end{equation}

Note that if $\hat G_n$ is a set estimator and $G$ is a measurable subset of $[0,1]^d$, then the quantity $|G\triangle\hat G_n|=\int_{[0,1]^d}|I(x\in\hat G_n)-I(x\in G)|dx$ is well defined and by Fubini's theorem it is measurable with respect to the probability measure $\PP_G$. Therefore one can measure the accuracy of the set estimator $\hat G_n$ on a given class of sets in the minimax framework : the risk of $\hat G_n$ on a class $\mathcal C$ is defined as 
\begin{equation*}
	\mathcal R_n(\hat G_n ; \mathcal C) = \sup_{G\in\mathcal C}\mathbb E_G[|G\triangle\hat G_n|].
\end{equation*}
For all the estimators that we will define in the sequel we will be interested in upper bounds on their risk, which give information about the rate at which these risks tend to zero, when the number $n$ of available observations tends to infinity. For a given class of subsets $\mathcal C$, the minimax risk on this class when $n$ observations are available is defined as 
\begin{equation*}
\mathcal R_n(\mathcal C)=\inf_{\hat G_n} \mathcal R_n(\hat G_n ; \mathcal C),
\end{equation*}
where the infimum is taken over all set estimators depending on $n$ observations. If $\mathcal R_n(\mathcal C)$ converges to zero, we call minimax rate of convergence on the class $\mathcal C$ the speed at which $\mathcal R_n(\mathcal C)$ tends to zero.

In this paper, we study minimax rates of convergence on two classes of subsets of $[0,1]^d$ : the class of all compact and convex sets, and the class of all polytopes with at most $r$ vertices, where $r$ is a given positive integer. 
Let $\mathcal C$ be a given class of subsets of $[0,1]^d$. We aim to provide with lower bounds on the minimax risks on the class $\mathcal C$. This lower bound can give much information on how close the risk of a given estimator is to the minimax risk on the class that we consider. If the rate (a sequence depending on $n$) of the upper bound on the risk of a given estimator matches with the rate of the lower bound on the minimax risk on the class $\mathcal C$, then this estimator is said to have the minimax rate of convergence on this class. 

We denote by $\rho$ the Euclidean distance in $\mathbb R^d$, by $B_d(y,r)$ the $d$-dimensional closed Euclidean ball centered at $y\in\mathbb R^d$ with radius $r$, and by $\beta_d$  the volume of the Euclidean unit ball in dimension $\mathbb R^d$. 
For any positive real number $x$, we denote by $\lfloor x \rfloor$ the greatest integer that is less or equal to $x$.
Any convex set that we will consider in the following is assumed to be compact and with nonempty interior in the considered topological space. 

\subsection{Former results and contributions}
Estimation of convex sets and, more generally, of sets, has been extensively studied in the previous decades (see the nice surveys given in Cuevas \cite{8'} and Cuevas and Fraiman \cite{8''} and the references therein, and related topics in \cite{18''0}). First works, in the 1960's, due to Renyi and Sulanke \cite{26}, \cite{27}, and Efron \cite{13} were motivated by issues of stochastic geometry, discussed, for instance, in the book by Kendall and Moran \cite{18'} and \cite{0}. Most of the works on estimation of convex sets dealt with models different than ours. Renyi and Sulanke, \cite{26}, \cite{27}, were the first to study the convex hull of a sample of $n$ i.i.d. random points in the plane. They obtained exact asymptotic formulas for the expected area and the expected number of vertices when the points are uniformly distributed over a convex set, and when they have a Gaussian distribution. They showed that if the points are uniformly distributed over a convex set $K$ in the plane $\mathbb R^2$, then the expected missing area $\mathbb E[|K\backslash\hat K|]$ of the convex hull $\hat K$ of the collection of these points is of the order 
\begin{itemize}
\item $n^{-2/3}$ if the boundary of $K$ is smooth,
\item $r\ln n /n$ if $K$ is a polygon with $r$ vertices.
\end{itemize}

This result was generalized to any dimension, and we refer to \cite{3} for an overview.

Estimation of convex sets in a multiplicative regression model has been investigated by Mammen and Tsybakov \cite{21'} and Korostelev and Tsybakov \cite{20}. The design $(X_1,\ldots,X_n)$ may be either random or deterministic, in $[0,1]^d$. In \cite{21'} Mammen and Tsybakov proposed an estimator of $G$ when it is assumed to be convex, based on likelihood-maximization over an $\varepsilon$-net, whose cardinality is bounded in terms of the metric entropy \cite{12}. They showed, with no assumption on the design, that the rate of their estimator cannot be improved.

The additive model \eqref{posprob} has been studied in \cite{19'} and \cite{20}, in the case where $G$ belongs to a smooth class of boundary fragments and the error terms are i.i.d. Gaussian variables with known variance. If $\gamma$ is the smoothness parameter of the studied class, it is shown that the rate of the minimax risk on the class is $n^{-\gamma/(\gamma+d-1)}$. The case of convex boundary fragments is covered by the case $\gamma=2$, which leads to the expected rate for the minimax risk, as we will discuss later (Section 5) : $n^{-2/(d+1)}$. It is important to note that in these works the authors always assumed that the fragment, which is included in $[0,1]^d$, has a boundary which is uniformly separated from $0$ and $1$. We will not make such an assumption in our work. Cuevas and Rodriguez-Cazal \cite{10}, and Pateiro Lopez \cite{24'}, studied the properties of set estimators of the support of a density, under several geometrical assumptions on the boundary of the unknown set.

One problem has not been investigated yet : how is the minimax rate of convergence modified if one assumes that the unknown set $G$, in model \eqref{posprob}, is a polytope, with a bounded number of vertices ? This question can be reformulated in a more general framework when one deals with boundary fragments : what is the minimax rate of convergence if $G$ is a fragment which belongs to a parametric family ? In the method used in \cite{19'} and \cite{20}, the true fragment is first approximated by an element of a parametric family of fragments, whose dimension is chosen afterwards according to the optimal bias-variance tradeoff, and the proposed estimator actually estimates the parametric approximation of the fragment $G$, and not directly $G$ itself. This idea is exploited in the present work, when we estimate convex sets, by using polytopal approximations. In the framework of fragments, the rate of convergence of the estimator when the target is the parametric fragment is found to be of the order $M/n$, where $M$ is the dimension of the parametric class of fragments. Again, the assumption of uniform separation from $0$ and $1$ is made. As we will show in the sequel, this assumption is essential in the parametric case, because if it is relaxed, an extra logarithmic factor appears in the rate.

In order to estimate convex sets, we will first approximate a convex set by a polytope, and then estimate that polytope. There is a wide literature on polytopal approximation of convex sets (cf. \cite{23}, \cite{14bis}, ...), which is of essential use in this paper.

For an integer $r\geq d+1$, we denote by $\mathcal P_r$ the class of all polytopes in $[0,1]^d$ with at most $r$ vertices. This class may be embedded into the finite-dimension space $\mathbb R^{dr}$ since any polytope is completely defined by the coordinates of its vertices. Therefore, one may expect that the problem of estimating $G\in\mathcal P_r$, for a given $r$, is parametric and therefore a rate of the order $1/n$ for the minimax risk $\mathcal R_n(\mathcal P_r)$, cf. \cite{18}. In Section 2, we propose an estimator that almost achieves this rate, up to a logarithmic factor. Moreover, we prove an exponential deviation inequality for the Nikodym distance between the estimator and the true polytope. Such an exponential inequality is of interest because it is much stronger than an upper bound on the risk of the estimator, and it is the key that leads to adaptive estimation, as we will see later. In Section 2, we show that this estimator has the minimax rate of convergence, so that the logarithmic factor in the rate is unavoidable. In Section 3, we extend the exponential deviation inequality of Section 2 in order to cover estimation of any convex set. 
In Section 4, we propose an estimator that is adaptive to the number of vertices of the estimated polytope, using as a convention that a non polytopal convex set can be considered as a polytope with infinitely many vertices. Section 5 is devoted to the proofs.

\section{Estimation of Convex Polytopes}

\subsection{Upper bound} 

We denote by $P_0$ the true polytope, i.e. $G=P_0$ in \eqref{posprob} and we assume that $P_0\in\mathcal P_r$. We denote by $\mathcal P_{r}^{(n)}$ the class of all the polytopes in $[0,1]^d$ with at most $r$ vertices with coordinates that are integer multiples of $\frac{1}{n}$. It is clear that the cardinality of $\mathcal P_{r}^{(n)}$ is less than $(n+1)^{dr}$. We have the following lemma, proved in Section 6.

\begin{lemma}
Let $r\leq n$. For any polytope $P$ in $\mathcal P_r$ there exists a polytope $P^*\in\mathcal P_{r}^{(n)}$ such that 
\begin{equation}
	\label{lem0}|P^*\triangle P|\leq \frac{2d^{d+1}(3/2)^d\beta_d}{n}.
\end{equation}
\end{lemma}

We estimate $P_0$ by a polytope in $\mathcal P_{r}^{(n)}$ that minimizes a given criterion. The criterion that we use is the sum of squared errors  
\begin{equation*}
\mathcal A(P,\{(X_i,Y_i)\}_{i=1,\ldots, n})=\sum_{i=1}^n(1-2Y_i)I(X_i\in P).
\end{equation*}
In order to simplify the notations, we will write $\mathcal A(P)$ instead of $\mathcal A(P,\{(X_i,Y_i)\}_{i=1,\ldots, n})$ in what follows. Note that if the noise terms $\xi_i, i=1,\ldots,n,$ are supposed to be Gaussian, then minimization of $\mathcal A(P)$ is equivalent to maximization of the likelihood.

Consider the set estimator of $P_0$ defined as  
\begin{equation}
	\label{defest} \hat P_n^{(r)}\in\underset{P\in\mathcal P_{r}^{(n)}}{\operatorname{argmin}} \text{ } \mathcal A(P).
\end{equation}
Note that since $\mathcal P_{r}^{(n)}$ is finite, the estimator $\hat P_n^{(r)}$ exists but is not necessarily unique.

Let us introduce the following assumption on the law of the noise terms $\xi_i, i=1,\ldots,n$ : 
\begin{AssumptionA}
The noise terms $\xi_i, i=1,\ldots,n,$ are subgaussian, i.e. satisfy the following exponential inequality :
\begin{equation*}
    \mathbb E[e^{u\xi_i}]\leq e^\frac{u^2\sigma^2}{2}, \forall u\in\mathbb R,
\end{equation*}
where $\sigma$ is a given positive number.
\end{AssumptionA}

Note that if the noise terms $\xi_i, i=1,\ldots,n,$ are i.i.d. centered Gaussian random variables, then Assumption A is satisfied.

The next theorem establishes an exponential deviation inequality for the estimator $\hat P_n^{(r)}$.
\begin{theorem}
Let Assumption A be satisfied.
For the estimator $\hat P_n^{(r)}$, there exist two positive constants $C_1$ and $C_2$, which depend on $d$ and $\sigma$ only, such that :  
\begin{equation*}
    \sup_{P\in\mathcal P_r}\mathbb P_{P}\Big[n\Big(|\hat P_n^{(r)}\triangle P|-\frac{2dr\ln n}{C_2n}\Big)\geq x\Big] \leq C_1e^{-C_2x}, \forall x>0.
\end{equation*}
\end{theorem}

The explicit forms of the constants $C_1$ and $C_2$ are given in the proof.
From the deviation inequality given in Theorem 1 one can easily derive that the risk of the estimator $\hat P_n^{(r)}$ on the class $\mathcal P_r$ is of the order $\frac{\ln n}{n}$. Indeed we have the following result.

\begin{corollary} Let the assumptions of Theorem 1 be satisfied. Then, for any positive number $q$, there exists a constant $A_q$ such that
\begin{equation*}
\sup_{P\in\mathcal P_r} \mathbb E_{P}\Big[|\hat P_n^{(r)}\triangle P|^q\Big]\leq A_q\left(\frac{dr\ln n}{n}\right)^q, \forall n\geq 1.
\end{equation*}
\end{corollary}
The explicit form of the constant $A_q$ can be easily derived from the proof.

\subsection{Lower bound} 
Corollary 1 gives an upper bound of the order $\frac{\ln n}{n}$ for the risk of our estimator $\hat P_n^{(r)}$. The next result shows that $\frac{\ln n}{n}$ is the minimax rate of convergence on the class $\mathcal P_r$.
\begin{theorem}
Assume that the noise terms $\xi_i, i=1,\ldots,n,$ are centered Gaussian random variables, with a given variance $\sigma^2>0$.
For every $r\geq d+1$, we have the following lower bound  
\begin{equation*}
      \label{LB}\inf_{\hat P}\sup_{P\in\mathcal P_r}\mathbb E_P\big[|\hat P\triangle P|\big] \geq \frac{\alpha^2\sigma^2\ln n}{n},
\end{equation*}
where $\displaystyle{\alpha=\frac{1}{2}-\frac{\ln 2}{2\ln 3}\approx 0.29...}$
\end{theorem}

Corollary 1 together with Theorem 3 gives the following bound on the class $\mathcal P_r$, in the case of Gaussian noise terms with variance $\sigma^2$ :
\begin{equation*}
	\alpha^2\sigma^2\leq\frac{n}{\ln n}\mathcal R_n(\mathcal P_r)\leq \frac{12dr}{1-e^{-\frac{1}{4\sigma^2}}}, 
\end{equation*}
for $n$ large enough and $d+1\leq r\leq n$. 
Note that the lower bound does not depend on the number of vertices $r$. This is because we prove our lower bound for the class $\mathcal P_{d+1}$ and we use that $\mathcal P_r\supseteq\mathcal P_{d+1}$, for $r\geq d+1$.

\section{Estimation of General Convex Sets}

\subsection{Upper bound}
Let us denote by $\mathcal C_d$ the class of all convex sets included in $[0,1]^d$.

Now we aim to estimate convex sets in the same model, without any assumption of the form of the unknown set. If $C$ is a convex set and $G=C$ in model \eqref{posprob}, an idea is to approximate $C$ by a convex polytope. For example one can select $r$ points on the boundary of $C$ and take their convex hull. This will give a polytope $C_r$ with $r$ vertices inscribed in $C$. In Section 2 we showed how to estimate such a $r$-vertex polytope as $C_r$. Thus, if $C_r$ approximates well $C$, an estimator of $C_r$ is a candidate to be a good estimator of $C$. The larger is $r$, the better $C_r$ should approximate $C$ with respect to the Nikodym distance defined in \eqref{defdist}. At the same time, when $r$ increases the upper bound given in Corollary 1 increases as well. Therefore $r$ should be chosen according to the bias-variance tradeoff.

For any integer $r\geq d+1$ consider again the estimator $\hat P_n^{(r)}$ defined in \eqref{defest}. However, now we chose a value for $r$ that depends on $n$ in order to achieve the bias-variance tradeoff.

\begin{theorem}
Consider model \eqref{posprob} with $G=C$, where $C$ is any convex subset of $[0,1]^d$.
Set $\displaystyle{r=\left\lfloor\left(\frac{n}{\ln n}\right)^\frac{d-1}{d+1}\right\rfloor}$, and let $\hat P_n^{(r)}$ the estimator defined in \eqref{defest}. 
Let Assumption A be satisfied.
Then, there exist positive constants $C_1, C_2$ and $C_3$, which depend on $d$ and $\sigma$ only, such that
\begin{equation*}
	 \sup_{C\in\mathcal C_d}\mathbb P_C\Big[n\Big(|\hat P_n^{(r)}\triangle C|-\left(\frac{C_3\ln n}{n}\right)^{2/(d+1)}\Big)\geq x\Big] \leq C_1e^{-C_2x}, \forall x>0.
\end{equation*}
\end{theorem}
The constants $C_1$ and $C_2$ are the same as in Theorem 1, and $C_3$ is given explicitly in the proof of the theorem. From Theorem 3 we get the next corollary.
\begin{corollary}
Let the assumptions of Theorem 3 be satisfied. Then, for any positive number $q$ there exists a positive constant $A'_q$ such that   
\begin{equation*}
\sup_{C\in\mathcal C_d} \mathbb E_C\Big[|\hat P_n^{(r)}\triangle C|^q\Big]\leq A'_q\left(\frac{\ln n}{n}\right)^\frac{2q}{d+1}, \forall n\geq 1.
\end{equation*}
\end{corollary} 
The explicit form of $A'_q$ can be easily derived from the proof.

\subsection{Lower bound}

In this section we give a lower bound on the minimax risk on the class $\mathcal C_d$ of all convex sets in $[0,1]^d$.  

\begin{theorem}
Assume that the noise terms $\xi_i, i=1,\ldots,n,$ are centered Gaussian random variables, with a given variance $\sigma^2>0$.
There exist a positive constant $C_{17}$ which depends only on the dimension $d$ and on $\sigma$, such that for any $n\geq 125$ and any estimator $\hat C$,
\begin{equation*}
	\sup_{C\in\mathcal C_d}\mathbb E_C\left[|C\triangle\hat C|\right] \geq C_{17}n^{-2/(d+1)}.
\end{equation*}
\end{theorem}
The explicit form of the constant $C_{17}$ can be found in the proof of the theorem. One can see that the lower bound given in Theorem 4 does not match the upper bound of in Theorem 3, where we had an extra logarithmic factor. Indeed we get that 
\begin{equation*}
	C_{17}n^{-2/(d+1)}\leq\mathcal R_n(\mathcal C_d)\leq 3\left(\frac{B_1\ln n}{n}\right)^\frac{2}{d+1}.
\end{equation*}

This gap is discussed in Section 5.

\section{Adaptive estimation}

In Section 2, we proposed an estimator that depends on the parameter $r$. A natural question is to find an estimator that is adaptive to $r$, i.e. that does not depend on $r$, but achieves the optimal rate on the class $\mathcal P_r$. 
The idea of the following comes from Lepski's method for adaptation (see \cite{21'0}, or \cite{8}, Section 1.5, for a nice overview).
Assume that the true number of vertices, denoted by $r^*$, is unknown, but is bounded from above by a given integer $R_n\geq d+1$ that may depend on $n$ and be arbitrarily large. Theorem 1 would provide the estimator $\hat P_n^{(R_n)}$, but it is clearly suboptimal if $r^*$ is small and $R_n$ is large. Indeed the rate of convergence of $\hat P_n^{(R_n)}$ is $\frac{R_n\ln n}{n}$, although the rate $\frac{r^*\ln n}{n}$ can be achieved according to Theorem 1, when $r^*$ is known. The procedure that we propose selects an integer $\hat r$ based on the observations, and the resulting estimator is $\hat P_n^{(\hat r)}$.

Note that $R_n$ should not be of order larger than $\left(\frac{n}{\ln n}\right)^{\frac{d-1}{d+1}}$, since for larger values of $r$, Corollaries 1 and 2 show that it is more efficient to consider the class $\mathcal C_d$ than the class $\mathcal P_r$.
Let us define : 
\begin{equation*}
	\hat r=\min\left\{ r\in\{d+1,\ldots,R_n\} : |\hat P_n^{(r)}\triangle\hat P_n^{(r')}|\leq \frac{6dr'\ln n}{C_2n}, \forall r'=r,\ldots,R_n\right\}.
\end{equation*}

The integer $\hat r$ is well defined, because the set in the brackets is not empty, since $R_n$ satisfies the condition.

Let us define the adaptive estimator $\hat P_n^{adapt}=\hat P_n^{(\hat r)}$. We then have the following theorem.


%
%
 
\begin{theorem}
Let Assumption A be satisfied.\\
Let $R_n=\lfloor \left(\frac{n}{\ln n} \right)^{\frac{d-1}{d+1}}\rfloor$ and $\phi_{n,r}=\min\left(\frac{r\ln n}{n},\left(\frac{\ln n}{n}\right)^{\frac{2}{d+1}}\right)$, for all integers $r\geq d+1$ and $r=\infty$.
There exists a positive constant $C_5$ that depends on $d$ and $\sigma$ only, such that the adaptive estimator $\hat P_n^{adapt}$ satisfies the following inequality :
\begin{equation*}
	\sup_{d+1\leq r\leq \infty}\sup_{P\in\mathcal P_r} \mathbb E_P\left[\phi_{n,r}^{-1}|\hat P_n^{adapt}\triangle P|\right] \leq C_5,
\end{equation*}
$\forall n\geq 1$, where $\mathcal P_{\infty}=\mathcal C_d$.
\end{theorem}

Thus, we show that one and the same estimator $\hat P_n^{adapt}$ attains the optimal rate simultaneously on all the classes $\mathcal P_r, d+1\leq r<\infty$, and near optimal rate (optimal up to a logarithmic factor) on the class $\mathcal C_d$ of all convex subsets of $[0,1]^d$. The explicit form of the constant $C_5$ can be easily derived from the proof of the theorem.


\section{Discussion}

In Theorems 3 and 4, the upper and lower bounds differ by a logarithmic factor, and a question is which of the two bounds could be improved. Theorems 1 and 2 show that the logarithmic factor is significant in the case of polytopes. Is it still the case for general convex sets ? 

Let us first answer the following question : what makes the estimation of sets on a given class $\mathcal C\subseteq\mathcal C_d$ difficult in the studied model ? First, it is the complexity of the class. As introduced by Dudley \cite{12}, the complexity of the class quantifies how big the class is, or in more precise words, the number of elements that are needed in order to discretize the class with a given precision. The more there are such elements, the more complex the class is, and the more complicated it is to estimate an unknown element of it. Second, it is how detectable the sets of the given class are, in our model. If the unknown subset $G$ is too small, then, with high probability, it contains no point of the design. Conditionally to this, all the data have the same distribution and no information in the sample can be used in order to detect $G$. A subset $G$ has to be large enough in order to be detectable by a given procedure. The threshold on the volume beyond which a subset cannot be detected by any procedure gives a lower bound on the rate of the minimax risk. In \cite{18'0}, Janson studied asymptotic properties on the maximal volume of holes with a given shape. A hole is a subset of $[0,1]^d$ that contains no point of the design $(X_1,\ldots,X_n)$. Janson showed that with high probability, there are convex and polytopal holes that have a volume of order $\ln n/n$. This result made it reasonable to think that $\ln n/n$ should be the order of a lower bound on the minimax risk in Theorem 2 ; this is the idea that we use in the proof of this theorem. The lower bound is attained on the polytopes with very small volumes. We do not use the specific structure of these polytopes to derive the lower bound : we only use the fact that some of them cannot be distinguished from the empty set, no matter what is the shape of their boundary, when we chose their volume of order no larger than $\frac{\ln n}{n}$. This shows that the rate $1/n$, which would come from the complexity of the parametric class $\mathcal P_r$, is not the right minimax rate of convergence ; the order $\ln n/n$, larger than $1/n$, imposes its law on this class.
On the other hand, the proof of our lower bound of the order $n^{-2/(d+1)}$ for general convex sets uses only the structure and regularity of the boundaries ; we do not deal especially with small hypotheses. The order $n^{-2/(d+1)}$ is much larger that $\ln n/n$, and therefore seems to determine the best lower bound achievable on the minimax risk on the class $\mathcal C_d$.

Regarding this discussion we formulate two conjectures.
\paragraph{Conjecture 1}
We conjecture that the risk of our estimator could be more sharply bounded than in Theorem 1, i.e. that
\begin{equation*}
	\max\left(\frac{\lambda_1\ln n}{n},\frac{\lambda_2r}{n}\right)\leq \mathcal R_n(\hat P_n^{(r)} ; \mathcal P_r) \leq \max\left(\frac{\lambda_3\ln n}{n},\frac{\lambda_4r}{n}\right),
\end{equation*}
for some positive constants $\lambda_1, \lambda_2, \lambda_3$ and $\lambda_4$. 
If $\ln n$ is sufficiently larger than $r$, the right order of the minimax risk is $\frac{\ln n}{n}$. If not, i.e. if the number of vertices of the unknown polytope can be large, the order of the risk is $\frac{r}{n}$. This lower bound is actually easy to prove when $d=2$, using the same scheme as in the proof of the case $d=2$ of Theorem 4.

\paragraph{Conjecture 2}
Let $\mu_0$ be a given positive number. If one considers the subclass $\mathcal P'_r(\mu_0)=\{P\in\mathcal P_r : |P|\geq\mu_0\}$, then subsets of $[0,1]^d$ with too small volume are not allowed anymore. Therefore, the hypotheses used in the proof of Theorem 4 are not valid anymore and we expect the minimax rate of convergence on this class to be of the order $1/n$.

\begin{remark}
If Conjecture 1 is true, and if we keep our method for estimating general convex sets and follow the proof of Theorem 3, the bias-variance tradeoff leads to a choice for $r$ of the order $n^{1/(d+1)}$, which is much larger than $\ln n$. Therefore, the risk has the rate $\frac{r}{n}=n^{-2/(d+1)}$ and the logarithmic factor is dropped.
\end{remark}

\section{Proofs}

\paragraph{Proof of Theorem 1} 
Let $P_0\in\mathcal P_r$ be the true polytope. Note that for all $\epsilon>0$, 
\begin{equation}
  \label{step1} \mathbb P_{P_0}\Big[|\hat P_n^{(r)}\triangle P_0|\geq \epsilon\Big] = \mathbb P_{P_0}\Big[\exists P\in\mathcal P_{r}^{(n)} : \mathcal A(P)\leq \mathcal A(P^*), |P\triangle P_0|\geq \epsilon\Big],
\end{equation}
where $P^*$ is a polytope chosen in $\mathcal P_{r}^{(n)}$ such that $|P^*\backslash P_0|\leq \frac{(4d)^{d+1}\beta_d}{n}$, cf. \eqref{lem0}. 
For any $P$ we have, by a simple algebra,
\begin{equation}
	\label{step1'''}\mathcal A(P)-\mathcal A(P^*) = \sum_{i=1}^nZ_i,
\end{equation}
where  
\begin{align*}
Z_i =  & I(X_i\in P)-I(X_i\in P^*)-2I(X_i\in P_0)\left[I(X_i\in P)-I(X_i\in P^*)\right] \nonumber \\
	 & -2\xi_i\left[I(X_i\in P)-I(X_i\in P^*)\right], \mbox{ } i=1,\ldots, n.
\end{align*}
The random variables $Z_i$ depend on $P$ but we omit this dependence in the notation. 
Therefore \eqref{step1} implies that 
\begin{align}
\mathbb P_{P_0}\Big[|\hat P_n^{(r)}\triangle P_0|\geq \epsilon\Big] & \leq \sum_{P\in \mathcal P_{r}^{(n)} :  |P\triangle P_0|\geq \epsilon}\mathbb P_{P_0}\Big[\sum_{i=1}^n Z_i\leq 0\Big] \nonumber \\
	& \label{step2} \leq  \sum_{P\in \mathcal P_{r}^{(n)} :  |P\triangle P_0|\geq \epsilon}\mathbb E_{P_0}\Big[\exp{(-u\sum_{i=1}^n Z_i)}\Big],
\end{align}
for all positive number $u$, by Markov's inequality. Since $Z_i$'s are mutually independent, we obtain 
\begin{equation}
  \label{step2'}
  \mathbb P_{P_0}\Big[|\hat P_n^{(r)}\triangle P_0|\geq \epsilon\Big] \leq  \sum_{P\in \mathcal P_{r}^{(n)} :  |P\triangle P_0|\geq \epsilon} \prod_{i=1}^n\mathbb E_{P_0}\Big[\exp{(-u Z_i)}\Big].
\end{equation}
By conditioning on $X_1$ and denoting by $W=I(X_1\in P)-I(X_1\in P^*)$ we have
\begin{align}
	\mathbb E_{P_0}\Big[ \exp(-uZ_1)\Big]  & = \mathbb E_{P_0}\Big[ \mathbb E_{P_0}\Big[ \exp(-uZ_1)|X_1\Big]\Big] \nonumber \\
	& = \mathbb E_{P_0}\Big[\exp\big( -uW +2uI(X_1\in P_0)W\big)\mathbb E_{P_0}\Big[ \exp\left(2u\xi_1W\right)|X_1\Big]\Big] \nonumber \\
	& = \mathbb E_{P_0}\Big[\exp\big( -uW+2uI(X_i\in P_0)W\big) \exp\big(2\sigma^2u^2I(X_1\in P\triangle P^*) \big)\Big] \nonumber \\
	\label{lastexpr} & = \mathbb E_{P_0}\Big[\exp\big( 2\sigma^2u^2I(X_1\in P\triangle P^*)-uW +2uI(X_1\in P_0)W\big)\Big].
\end{align}
We will now reduce the last expression in \eqref{lastexpr}. It is convenient to use the following table in which the first three columns represent the values that can be taken by the binary variables $I(X_1\in P)$, $I(X_1\in P^*)$ and $I(X_1\in P_0)$ respectively, and the last column gives the resulting value of the term $\exp\big( 2\sigma^2u^2I(X_1\in P\triangle P^*)-uW+2uI(X_1\in P_0)W\big)$ that is under the expectation in \eqref{lastexpr}.

{\renewcommand{\arraystretch}{1.5}
\begin{table}[h]
\centering
\begin{tabular}{|c|c|c|c|}
  \hline
 $P$ & $P^*$ & $P_0$ & Value \\
  \hline\hline
  $1$ &  $1$ & $1$ & $1$ \\
  \hline
  $1$ &  $1$ & $0$ & $1$ \\
  \hline
  $1$ &  $0$ & $1$ & $\exp(2\sigma^2u^2+u)$ \\
  \hline
  $1$ &  $0$ & $0$ & $\exp(2\sigma^2u^2-u)$ \\
  \hline
  $0$ &  $1$ & $1$ & $\exp(2\sigma^2u^2-u)$ \\
  \hline
  $0$ &  $1$ & $0$ & $\exp(2\sigma^2u^2+u)$ \\
  \hline
  $0$ &  $0$ & $1$ & $1$ \\
  \hline
  $0$ &  $0$ & $0$ & $1$ \\
  \hline
\end{tabular}
\end{table}}

Hence one can write
\begin{align*}
	\mathbb E_{P_0}\Big[ \exp(-uZ_1)\Big] = & \mbox{ }1-|P\triangle P^*|+e^{2\sigma^2u^2+u}\big(\left|(P\cap P_0)\backslash P^* \right| +\left|P^*\backslash (P\cup P_0) \right| \big) \\
	&+e^{2\sigma^2u^2-u}\big(\left|(P^*\cap P_0)\backslash P \right|+\left|P\backslash (P^*\cup P_0) \right| \big).
\end{align*}
Besides by the triangle inequality,
\begin{equation*}
	|P\triangle P_0| \leq |P\triangle P^*|+|P^*\triangle P_0|,
\end{equation*}
which implies
\begin{align}
	\mathbb E_{P_0}\Big[ \exp(-uZ_1)\Big] \leq & \mbox{ }1-|P\triangle P_0|+ |P^*\triangle P_0|    +e^{2\sigma^2u^2+u}\big(\left|P_0\backslash P^* \right| +\left|P^*\backslash P_0 \right| \big) \nonumber \\
	&+e^{2\sigma^2u^2-u}\big(\left|P_0\backslash P \right|+\left|P\backslash P_0 \right| \big) \nonumber \\
	\label{step1''''}\leq & \mbox{ } 1-|P\triangle P_0|+|P^*\triangle P_0| +e^{2\sigma^2u^2+u}|P^*\triangle P_0|+e^{2\sigma^2u^2-u}|P\triangle P_0| \\
	\leq & \mbox{ } 1-|P\triangle P_0|\left(1-e^{2\sigma^2u^2-u}\right)+\frac{2d^{d+1}(3/2)^d\beta_d}{n}\left(1+e^{2\sigma^2u^2+u}\right). \nonumber
\end{align}
Choose $u=\frac{1}{4\sigma^2}$. Then the quantity $1-e^{2\sigma^2u^2-u}$ is positive and if $|P\triangle P_0|\geq\epsilon$, then
\begin{equation}
	\label{lasteqq}\mathbb E_{P_0}\Big[ \exp(-uZ_1)\Big] \leq 1-\epsilon\left(1-e^{-\frac{1}{4\sigma^2}}\right)+\frac{2d^{d+1}(3/2)^d\beta_d}{n}\left(1+e^{\frac{3}{8\sigma^2}}\right).
\end{equation}
We set $\tilde{C_1}=1+e^{\frac{3}{8\sigma^2}}$ and $C_2=1-e^{-\frac{1}{4\sigma^2}}$. These are positive constants that do not depend on $n$ or $P_0$. From \eqref{step2'} and \eqref{lasteqq}, and by the independence of $Z_i$'s we have
\begin{align}
	\mathbb P_{P_0}\Big[|\hat P_n^{(r)}\triangle P_0|\geq \epsilon\Big] & \leq \sum_{P\in \mathcal P_{r}^{(n)} :  |P\triangle P_0|\geq \epsilon} \Big(1-C_6\epsilon+\frac{2d^{d+1}(3/2)^d\beta_d\tilde{C_1}}{n}\Big)^n \nonumber \\
	& \leq (n+1)^{dr}\Big(1-C_2\epsilon+\frac{2d^{d+1}(3/2)^d\beta_d\tilde{C_1}}{n}\Big)^n \nonumber \\ 
	& \leq \exp\left(dr\ln (n+1)-C_2\epsilon n+2d^{d+1}(3/2)^d\beta_d\tilde{C_1}\right) \nonumber \\
	& \label{31415} \leq  \exp\left(2dr\ln n-C_2\epsilon n+C_8\right),
\end{align}
where $C_1=\exp\left(2d^{d+1}(3/2)^d\beta_d\tilde{C_1}\right)$, noting that $n+1\leq n^2$.
Therefore if we set $\epsilon=\frac{2dr\ln n}{C_2n}+\frac{x}{n}$ for a positive number $x$, we get the following deviation inequality
\begin{equation*}
\mathbb P_{P_0}\left[n\left(|\hat P_n^{(r)}\triangle P_0|-\frac{2dr\ln n}{C_2n}\right)\geq x\right]\leq C_1e^{-C_2x}.
\end{equation*}
$\blacksquare$

\paragraph{Proof of Corollary 1} 
Corollary 1 follows directly from Theorem 1 and Fubini's theorem. Indeed, if we denote by $Z:=|\hat P_n^{(r)}\triangle P_0|$ and by $\mathbb P_Z$ its distribution measure, then $Z$ is a continuous and nonnegative random variable and we have, by Fubini's theorem, that
\begin{align*}
	\mathbb E_{P_0}[Z^q] & = q\int_0^\infty u^{q-1}\mathbb P_Z[Z\geq u]du \nonumber  \\
	& \leq q\int_0^\frac{2dr\ln n}{C_2n}u^{q-1}du + q\int_0^\infty \left(u+\frac{2dr\ln n}{C_2n}\right)^{q-1}\mathbb P_Z\left[Z\geq u+\frac{2dr\ln n}{C_2n}\right]du \nonumber \\
	& = \left(\frac{2dr\ln n}{C_2n}\right)^q + q\int_0^\infty \left(u+\frac{2dr\ln n}{C_2n}\right)^{q-1}\mathbb P_Z\left[n\left(Z-\frac{2dr\ln n}{C_2n}\right)\geq nu\right]du \nonumber \\
	& \leq \left(\frac{2dr\ln n}{C_2n}\right)^q + q\int_0^\infty \left(u+\frac{2dr\ln n}{C_2n}\right)^{q-1}C_1e^{-C_2nu}du \mbox{ }\mbox{ by Theorem 1} \nonumber \\ 
	& \leq \left(\frac{2dr\ln n}{C_2n}\right)^q + C_1q\max(1,2^{q-1})\int_0^\infty \left(v^{q-1}+\left(\frac{2dr\ln n}{C_2n}\right)^{q-1}\right)e^{-C_2nv}dv \nonumber \\
	& \leq 3\left(\frac{2dr\ln n}{C_2n}\right)^q,
\end{align*}
for $n$ large enough. Note that the sixth step of this proof comes from the easy fact that for any positive numbers $a$ and $b$, $(a+b)^{q-1}\leq 2^{q-1}(a^{q-1}+b^{q-1})$ if $q-1>0$, and $(a+b)^{q-1}\leq a^{q-1}+b^{q-1}$ if $q-1\leq 0$, and the seventh comes from the fact that $\int_0^\infty v^{q-1}e^{-q}dv = (q-1)!$.
$\blacksquare$

\paragraph{Proof of Theorem 2}
This proof is a simple application of Corollary 2.6 in \cite{30}. Let $M$ be a positive integer, and $h=\frac{1}{M+1}$. Let $T_k, k=1,\ldots, M$ be $M$ disjoint polytopes in $\mathcal P_{d+1}$ and with same volume : $|T_1|=\ldots=|T_M|=h/2$, where $h=M^{-1}$.

For $k=1,\ldots,M$ we denote by $\mathbb P_k$ the probability distribution of the observations $(X_i,Y_i), i=1,\ldots, n$ when $G=T_k$ in \eqref{posprob}, and by $\mathbb E_k$ the expectation with respect to this distribution. A simple computation shows that the Kullback-Leibler divergence $K(\mathbb P_k,\mathbb P_l)$ between $\mathbb P_k$ and $\mathbb P_l$, for $k\ne l$, is equal to $\frac{nh}{4\sigma^2}$. On the other hand, the distance between $T_k$ and $T_l$, for $k\ne l$, is $|T_k\triangle T_l|=|T_k|+|T_l|=h$. Then
\begin{equation*}
	\frac{1}{M+1}\sum_{j=1}^MK(\mathbb P_j,\mathbb P_0) = \frac{Mnh}{4(M+1)\sigma^2} \leq \frac{n}{4M\sigma^2}.
\end{equation*}
Let $\alpha\in(0,1)$, and $\gamma=\frac{1}{2\sigma^2\alpha}$. Then, if $M=\frac{\gamma n}{\ln n}$, supposed without loss of generality to be an integer, we have 
\begin{equation*}
	4\sigma^2\alpha M\ln M = 2n+2\ln \gamma\frac{\ln n}{n}-2n\frac{\ln\ln n}{n} \geq n
\end{equation*}
for $n$ large enough, so that
\begin{equation*}
	\frac{1}{M+1}\sum_{j=1}^MK(\mathbb P_j,\mathbb P_0) \leq \alpha\ln M.
\end{equation*}
Therefore, applying Corollary 2.6 in \cite{30} with the pseudo distance defined in \eqref{defdist}, we set for $r\geq d+1$ the following inequality 
\begin{equation*}
\inf_{\hat P}\sup_{P\in\mathcal P_{d+1}}\mathbb E_P\big[|\hat P\triangle P|\big] \geq \frac{1}{M+1}\Big(\frac{\ln{(M+1)}-\ln 2}{\ln M} -\alpha\Big).
\end{equation*}
For $n$ great enough we have $M\geq3$ and $\frac{\ln{(M+1)}-\ln 2}{\ln M}\geq 1-\frac{\ln2}{\ln3}$. We choose $\alpha = \frac{1}{2}-\frac{\ln2}{2\ln3}\in(0,1)$. So, we get
\begin{equation*}
    \inf_{\hat P}\sup_{P\in\mathcal P_{d+1}}\mathbb E_P\big[|\hat P\triangle P|\big] \geq \frac{\alpha}{M+1} \geq \frac{\alpha}{2M} \geq \frac{\alpha\ln n}{\gamma n} \geq \frac{\alpha^2\sigma^2\ln n}{n}.
\end{equation*}
This immediately implies Theorem 2.
$\blacksquare$

\paragraph{Proof of Theorem 3}

The idea of the proof is very similar to that of Theorem 1. Here we need to control an extra bias term, due to the approximation of $C$ by a $r$-vertex polytope. We give the following lemma (cf. \cite{14bis}).

\begin{lemma} 
Let $r\geq d+1$ be a positive integer. For any convex set $C\subseteq \mathbb R^d$ there exists a polytope $C_r$ with at most $r$ vertices such that
\begin{equation*}
	|C\triangle C_r|\leq Ad\frac{|C|}{r^{2/(d-1)}},
\end{equation*}
where $A$ is a positive constant that does not depend on $r, d$ and $C$.
\end{lemma}
Let $P^*$ be a polytope chosen in $\mathcal P_{r}^{(n)}$ such that$|P^*\triangle C_r|\leq \frac{(4d)^{d+1}\beta_d}{n}$, like in the proof of Theorem 1. Thus by the triangle inequality, 
\begin{equation*}
	|P^*\triangle C| \leq |P^*\triangle C_r|+|C_r\triangle C| \leq \frac{Ad}{r^{2/(d-1)}}+\frac{(4d)^{d+1}\beta_d}{n}.
\end{equation*}
We now bound from above the probability $\displaystyle{\mathbb P_C\Big[|\hat P_n^{(r)}\triangle C|\geq \epsilon\Big]}$ for any $\epsilon>0$. 
As in \eqref{step1} and \eqref{step2} we have  
\begin{align*}
	\mathbb P_C\Big[|\hat P_n^{(r)}\triangle C|\geq \epsilon\Big] & \leq \mathbb P_C\Big[\exists P\in\mathcal P_{2r}^{(n)}, \mathcal A(P)\leq \mathcal A(P^*), |P\triangle C|\geq \epsilon\Big] \\
	& \leq  \sum_{P\in \mathcal P_{2r}^{(n)} :  |P\triangle C|\geq \epsilon} \mathbb P_C\Big[\mathcal A(P)\leq \mathcal A(P^*)\Big].	
\end{align*}
Repeating the argument in \eqref{step1'''} with $C$ instead of $P_0$ we set
\begin{equation*}
	\mathcal A(P)-\mathcal A(P^*) = \sum_{i=1}^nZ_i,
\end{equation*}
where
\begin{align*}
	Z_i =  & I(X_i\in P)-I(X_i\in P^*)-2I(X_i\in C)\left[I(X_i\in P)-I(X_i\in P^*)\right] \\
	 & -2\xi_i\left[I(X_i\in P)-I(X_i\in P^*)\right], \mbox{ } i=1,\ldots, n.
\end{align*}
The rest of the proof is very similar to the one of Theorem 1. Indeed, replacing $P_0$ by $C$ in that proof between \eqref{step1} and \eqref{step1''''}, and $\frac{2d^{d+1}(3/2)^d\beta_d}{n}$ by $\frac{2d^{d+1}(3/2)^d\beta_d}{n}+\frac{Ad}{r^{2/(d-1)}}$ in \eqref{lasteqq} and \eqref{31415} one gets :

\begin{align*}
	\mathbb P_C\Big[|\hat P_n^{(r)}\triangle C|\geq \epsilon\Big] & \leq \sum_{P\in \mathcal P_{r}^{(n)} :  |P\triangle C|\geq \epsilon} \left(1-C_2\epsilon+\tilde{C_1}\left(\frac{Ad}{r^{2/(d-1)}}+\frac{2d^{d+1}(3/2)^d\beta_d}{n}\right)\right)^n \\
	& \leq (n+1)^{dr}\left(1-C_2\epsilon+\tilde{C_1}\left(\frac{Ad}{r^{2/(d-1)}}+\frac{2d^{d+1}(3/2)^d\beta_d}{n}\right)\right)^n \\ 
	& \leq \exp\left(2dr\ln n-C_2\epsilon n+\tilde{C_1}\left(\frac{Adn}{r^{2/(d-1)}}+2d^{d+1}(3/2)^d\beta_d\right) \right).
\end{align*}
Therefore if we set $\epsilon=\frac{2dr\ln n}{C_2n}+\frac{\tilde{C_1}Ad}{C_2r^{2/(d-1)}}+\frac{x}{n}$ for a positive number $x$, we get the following deviation inequality
\begin{equation*}
\mathbb P_C\left[n\left(|\hat P_n^{(r)}\triangle C|-\frac{2dr\ln n}{C_2n}-\frac{\tilde{C_1}Ad}{C_2r^{2/(d-1)}}\right)\geq x\right]\leq C_1e^{-C_2x},
\end{equation*}
where the constants are defined as in the previous section.
That ends the proof of Theorem 3 by choosing $r=\lfloor\left(\frac{n}{\ln n}\right)^\frac{d-1}{d+1}\rfloor$, and the constant $C_3$ is given by 
\begin{equation*}
	C_3 = \left(1+\tilde{C_1}A\right)\frac{d}{C_2}=\left(1+(1+e^{3/(8\sigma^2)})A\right)\frac{d}{1-e^{-1/(4\sigma^2)}}. \mbox{ }\mbox{ }\mbox{ }\mbox{ }\mbox{ }\mbox{ }\mbox{ }\mbox{ }\mbox{ }\mbox{ }\mbox{ }\mbox{ }\mbox{ }\mbox{ }\mbox{ }\mbox{ }\mbox{ }\mbox{ }\mbox{ }\mbox{ }\mbox{ }\mbox{ }\mbox{ }\mbox{ }\mbox{ }\mbox{ }\mbox{ }\mbox{ }\mbox{ }\mbox{ }\mbox{ }\mbox{ }\mbox{ }\mbox{ }\mbox{ }\mbox{ }\mbox{ }\mbox{ }\blacksquare
\end{equation*}

\paragraph{Proof of Theorem 4}
We first prove this theorem in the case $d=2$ and then generalize the proof for $d\geq 3$.\\
We more or less follow the lines of the proof of the lower bound in \cite{21} (which is similar to the proof of Assouad's lemma, see \cite{30}). Let $G$ be the disk centered in $(1/2,1/2)$ of radius $1/2$, and $P$ be a regular convex polygon with $M$ vertices, all of them lying on the edge of $G$. Each edge of $P$ cuts a cap off $G$, of area $h$, with $\pi^3/(12M^3)\leq h\leq \pi^3/M^3$ as soon as $M\geq 6$, which we will assume in the sequel. We denote these caps by $D_1, \ldots, D_M$, and for any $\omega=(\omega_1,\ldots,\omega_M)\in\{0,1\}^M$ we denote by $G_\omega$ the set made of $G$ out of which we took all the caps $D_j$ for which $\omega_j=0$, $j=1,\ldots,M$. 

For $j=1,\ldots,M$, and $(\omega_1,\ldots,\omega_{j-1},\omega_{j+1},\ldots,\omega_M)\in\{0,1\}^{M-1}$ we denote by  
\begin{equation*}
\omega^{(j,0)}=(\omega_1,\ldots,\omega_{j-1},0,\omega_{j+1},\ldots,\omega_M) \mbox{  and by}
\end{equation*}
\begin{equation*}
\omega^{(j,1)}=(\omega_1,\ldots,\omega_{j-1},1,\omega_{j+1},\ldots,\omega_M).
\end{equation*} 
Therefore note that for any $j=1,\ldots,M$, and $(\omega_1,\ldots,\omega_{j-1},\omega_{j+1},\ldots,\omega_M)\in\{0,1\}^{M-1}$, 
\begin{equation*}
	|G_{\omega^{(j,0)}}\triangle G_{\omega^{(j,1)}}|=h.
\end{equation*}

For two probability measures $\mathbb P$ and $\mathbb Q$ defined on the same probability space and having densities denoted respectively by $p$ and $q$ with respect to a common measure $\nu$ (we also denote by $d\mathbb P = pd\nu$ and $d\mathbb Q = qd\nu$), we call $H(\mathbb P,\mathbb Q)$ the Hellinger distance between $\mathbb P$ and $\mathbb Q$, defined as 
\begin{equation*}
H(\mathbb P,\mathbb Q)=\left(\int(\sqrt p-\sqrt q)^2\right)^{1/2}.
\end{equation*}
Some useful properties of the Hellinger distance can be found in \cite{30}, Section 2.4.

Now, let us consider any estimator $\hat G$. For $j=1,\ldots,M$ we denote by $A_j$ the smallest convex cone with origin at $(1/2,1/2)$ and which contains the cap $D_j$. Note that the cones $A_j, j=1,\ldots,M$ have pairwise a null Lebesgue measure intersection. Then, we have the following inequalities :

\begin{align}
	& \sup_{G\in\mathcal C_2}\mathbb E_G\left[|G\triangle \hat G|\right] \nonumber \\ 
	& \geq \frac{1}{2^M}\sum_{\omega\in\{0,1\}^M}\mathbb E_{G_\omega}\left[|G_\omega\triangle \hat G|\right] \nonumber \\
	& \geq \frac{1}{2^M}\sum_{\omega\in\{0,1\}^M}\sum_{j=1}^M\mathbb E_{G_\omega}\left[|(G_\omega\cap A_j)\triangle (\hat G\cap A_j)|\right] \nonumber \\
	& = \frac{1}{2^M}\sum_{j=1}^M\sum_{\omega\in\{0,1\}^M}\mathbb E_{G_\omega}\left[|(G_\omega\cap A_j)\triangle (\hat G\cap A_j)|\right] \nonumber 
\end{align}
\begin{align}
	& = \frac{1}{2^M}\sum_{j=1}^M\underset{\omega_1,\ldots,\omega_{j-1},\omega_{j+1},\ldots,\omega_M}{\sum\ldots\sum}\Big(\mathbb E_{G_\omega^{(j,0)}}\left[|(G_\omega^{(j,0)}\cap A_j)\triangle (\hat G\cap A_j)|\right] \nonumber \\
	\label{25121986} & \mbox{ }\mbox{ }\mbox{ }\mbox{ }\mbox{ }\mbox{ }\mbox{ }\mbox{ }\mbox{ }\mbox{ }\mbox{ }\mbox{ }\mbox{ }\mbox{ }\mbox{ }\mbox{ }\mbox{ }\mbox{ }\mbox{ }\mbox{ }\mbox{ }\mbox{ }\mbox{ }\mbox{ }\mbox{ }\mbox{ }\mbox{ }\mbox{ }\mbox{ }\mbox{ }\mbox{ }\mbox{ }\mbox{ }\mbox{ }\mbox{ }\mbox{ }\mbox{ }\mbox{ }\mbox{ }\mbox{ }+ \mathbb E_{G_\omega^{(j,1)}}\left[|(G_\omega^{(j,1)}\cap A_j)\triangle (\hat G\cap A_j)|\right]\Big). 
\end{align}
Besides for any $j=1,\ldots,M$ and $(\omega_1,\ldots,\omega_{j-1},\omega_{j+1},\ldots,\omega_M)\in\{0,1\}^{M-1}$ we have 

\begin{align}
	& \mathbb E_{G_\omega^{(j,0)}}\left[|(G_\omega^{(j,0)}\cap A_j)\triangle (\hat G\cap A_j)|\right]+\mathbb E_{G_\omega^{(j,1)}}\left[|(G_\omega^{(j,1)}\cap A_j)\triangle (\hat G\cap A_j)|\right] \nonumber \\
	& = \int_{\left([0,1]^2\times\mathbb R\right)^{n}}|(G_\omega^{(j,0)}\cap A_j)\triangle (\hat G\cap A_j)|d\mathbb P_{G_\omega^{(j,0)}}^{\otimes n} \nonumber \\
	& \mbox{}\mbox{ }\mbox{ }\mbox{ }\mbox{ }\mbox{ }\mbox{ }\mbox{ }\mbox{ }\mbox{ }\mbox{ }\mbox{ }\mbox{ }\mbox{ }\mbox{ }\mbox{ }\mbox{ }\mbox{ }\mbox{ }\mbox{ }\mbox{ }\mbox{ }\mbox{ }+\int_{\left([0,1]^2\times\mathbb R\right)^{n}}|(G_\omega^{(j,1)}\cap A_j)\triangle (\hat G\cap A_j)|d\mathbb P_{G_\omega^{(j,1)}}^{\otimes n} \nonumber \\
\end{align}
\begin{align}
	& \geq \int_{\left([0,1]^2\times\mathbb R\right)^{n}}\left(|(G_\omega^{(j,0)}\cap A_j)\triangle (\hat G\cap A_j)| + |(G_\omega^{(j,1)}\cap A_j)\triangle (\hat G\cap A_j)|\right) \min(d\mathbb P_{G_\omega^{(j,0)}}^{\otimes n},d\mathbb P_{G_\omega^{(j,1)}}^{\otimes n}) \nonumber \\
	& \geq \int_{\left([0,1]^2\times\mathbb R\right)^{n}}\left(|(G_\omega^{(j,0)}\cap A_j)\triangle (G_\omega^{(j,1)}\cap A_j)|\right)\min(d\mathbb P_{G_\omega^{(j,0)}}^{\otimes n},d\mathbb P_{G_\omega^{(j,1)}}^{\otimes n}), \nonumber \\
	&  \mbox{ }\mbox{ }\mbox{ }\mbox{ }\mbox{ }\mbox{ }\mbox{ }\mbox{ }\mbox{ }\mbox{ }\mbox{ }\mbox{ }\mbox{ }\mbox{ }\mbox{ }\mbox{ }\mbox{ }\mbox{ }\mbox{ }\mbox{ }\mbox{ }\mbox{ }\mbox{ }\mbox{ }\mbox{ }\mbox{ }\mbox{ }\mbox{ }\mbox{ }\mbox{ }\mbox{ }\mbox{ }\mbox{ }\mbox{ }\mbox{ }\mbox{ }\mbox{ }\mbox{ }\mbox{ }\mbox{ }\mbox{ }\mbox{ }\mbox{ }\mbox{ }\mbox{ }\mbox{ }\mbox{ }\mbox{ }\mbox{ }\mbox{ }\mbox{ }\mbox{ }\mbox{ }\mbox{ }\mbox{ }\mbox{ }\mbox{ }\mbox{ }\mbox{ }\mbox{ }\mbox{ }\mbox{ }\mbox{ }\mbox{ }\mbox{ }\mbox{ }\mbox{ }\mbox{ }\mbox{ }\mbox{ }\mbox{ }\mbox{ }\mbox{ }\mbox{ }\mbox{ by the triangle inequality}\nonumber \\
	& = h\int_{\left([0,1]^2\times\mathbb R\right)^{n}}\min(d\mathbb P_{G_\omega^{(j,0)}}^{\otimes n},d\mathbb P_{G_\omega^{(j,1)}}^{\otimes n}) \nonumber \\
	& \geq \frac{h}{2}\left(1-\frac{H^2(\mathbb P_{G_\omega^{(j,0)}}^{\otimes n},\mathbb P_{G_\omega^{(j,1)}}^{\otimes n})}{2}\right)^2  \nonumber \\
	\label{25061986} & = \frac{h}{2}\left(1-\frac{H^2(\mathbb P_{G_\omega^{(j,0)}},\mathbb P_{G_\omega^{(j,1)}})}{2}\right)^{2n},
\end{align}
using properties of the Hellinger distance (cf. Section 2.4. in \cite{30}).
To compute the Hellinger distance between $\mathbb P_{G_\omega^{(j,0)}}$ and $\mathbb P_{G_\omega^{(j,1)}}$ we use the following lemma.
\begin{lemma}
For any integer $d\geq 2$, if $G_1$ and $G_2$ are two subsets of $[0,1]^d$, then 
\begin{equation*}
	H^2(\mathbb P_{G_1},\mathbb P_{G_2}) = 2(1-e^{-\frac{1}{8\sigma^2}})|G_0\triangle G_1|.
\end{equation*}
\end{lemma}
Then if we denote by $C_9=1-e^{-\frac{1}{8\sigma^2}}$, it follows from \eqref{25121986} and \eqref{25061986} that  
\begin{align*}
	\sup_{G\in\mathcal C_2}\mathbb E_G\left[|G\triangle \hat G|\right] & \geq \frac{1}{2^M}.M.2^{M-1}.\frac{h}{2}(1-C_9h)^{2n} \\
	& \geq \frac{Mh}{4}(1-C_9h)^{2n} \\
	& \geq \frac{\pi^3}{12M^2}(1-\pi^3C_9/M^3)^{2n}.
\end{align*}
Besides, since we assumed that $M\geq 6$, we have that 
\begin{equation*}
	\pi^3C_9/M^3\leq \pi^3C_9/6^3 =\frac{\pi^3}{6^3}\left(1-\exp(-\frac{1}{8\sigma^2})\right) \leq \frac{\pi^3}{6^3}<1, 
\end{equation*}
and we get by concavity of the logarithm  
\begin{equation*}
	\sup_{G\in\mathcal C_2}\mathbb E_G\left[|G\triangle \hat G|\right] \geq \frac{\pi^3}{12M^2}\exp\left(\frac{432\ln(1-\pi^3/216)\left(1-\exp(-\frac{1}{8\sigma^2})\right)nM^{-3}}{\pi^3}\right) \geq C_{14}n^{-2/3},
\end{equation*}
if we take $M=\lfloor n^{1/3} \rfloor$, where $\displaystyle{C_{14}=\frac{\pi^3}{12}\exp\left(\frac{432\ln(1-\pi^3/216)\left(1-\exp(-\frac{1}{8\sigma^2}\right)}{\pi^3}\right)}$ is a positive constant that depends only on $\sigma$. This inequality holds for $n\geq 216$, so that $M\geq 6$. \\

We now deal with the case $d\geq 3$. Let us first recall some definitions and resulting properties, that can also be found in \cite{19'0}.
\begin{definition} 
	Let $(S,\rho)$ be a metric space and $\eta$ a positive number. \\
	A family $\mathcal Y \subseteq S$ is called an $\eta$-packing family if and only if $\rho(y,y')\geq\eta$, for $(y, y')\in\mathcal Y$ with $y\neq y'$. \\
	An $\eta$-packing family is called maximal if and only if it is not strictly included in any other $\eta$-packing family . 
	A family $\mathcal Z$ is called an $\eta$-net if and only if for all $x\in S$, there is an element $z\in\mathcal Z$ which satisfies $\rho(x,z)\leq\eta$.
\end{definition}

We now give a Lemma.

\begin{lemma}
Let $S$ be the sphere with center $a_0=(1/2,\ldots,1/2)\in\mathbb R^d$ and radius $1/2$, and $\rho$ the Euclidean distance in $\mathbb R^d$. We still denote by $\rho$ its restriction on $S$. 

Let $\eta\in (0,1)$. Then any $\eta$-packing family of $(S,\rho)$ is finite, and any maximal $\eta$-packing family has a cardinality $M_\eta$ that satisfies the inequalities
\begin{equation}
	\label{Lemma6} \frac{d\sqrt{2\pi}}{2^{d-1}\sqrt{d+2}\eta^{d-1}}\leq M_\eta \leq \frac{4^{d-2}\sqrt{2\pi d}}{3^{(d-3)/2}\eta^{d-1}}.
\end{equation}
\end{lemma}

The construction of the hypotheses used for the lower bound in the case $d=2$ requires a little more work in the general dimension case, since it is not always possible to construct a regular polytope with a given number of vertices or facets, and inscribed in a given ball. For the following geometrical construction, we refer to Figure 1. 

Let $G_0$ be the closed ball in $\mathbb R^d$, with center $a_0=(1/2,\ldots,1/2)$ and radius $1/2$, so that $G_0\subseteq [0,1]^d$. Let $\eta\in (0,1)$ which will be chosen precisely later, and $\{y_1,\ldots,y_{M_\eta}\}$ a maximal $\eta$-packing family of $S=\partial G_0$. The integer $M_\eta$ satisfies \eqref{Lemma6} by Lemma 4.
For $j\in\{1,\ldots,M_\eta\}$, we set by $U_j=S\cap B_{d}(y_j,\eta/2)$, and denote by $W_j$ the $d-2$ dimensional sphere $S\cap \partial B_{d}(y_j,\eta/2)$. Let $H_j$ be affine hull of $W_j$, i.e. its supporting hyperplane. $H_j$ dissects the space $\mathbb R^d$ into two halfspaces. Let $H_j^-$ be the one that contains the point $y_j$. For $\omega=(\omega_1,\ldots,\omega_{M_\eta})\in\{0,1\}^{M_\eta}$, we denote by
\begin{equation*}
	G_\omega = G_0 \backslash (\bigcap_{j=1,\ldots,M_\eta : \omega_j=0}H_j^-).
\end{equation*}
The set $G_\omega$ is made of $G_0$ from which we remove all the caps cut off by the hyperplanes $H_j$, for all the indices $j$ such that $\omega_j=0$.

\begin{figure}
\centering
\includegraphics[width=90mm]{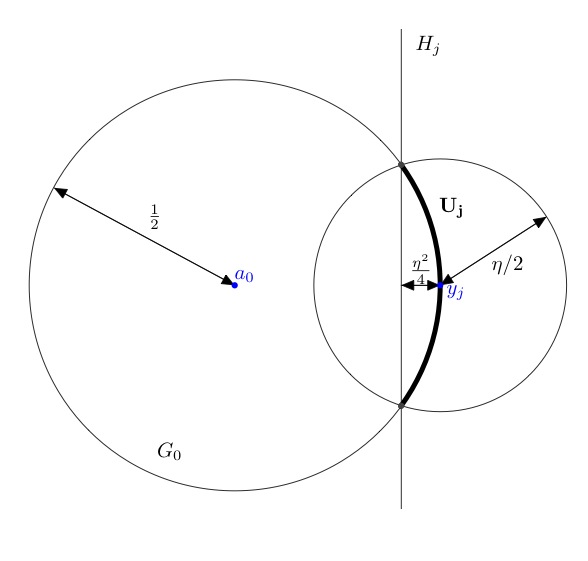}
\caption{Construction of the hypotheses}
\label{Fig 1}
\end{figure}

For each $j\in\{1,\ldots,M_{\eta}\}$, let $A_j$ be the smallest closed convex cone with vertex $a_0=(1/2,\ldots,1/2)$ that contains $U_j$. Note that the cones $A_j, j=1,\ldots,M_\eta$ have pairwise empty intersection, since $G_0$ is convex and the sets $U_j$ are disjoint. We are now all set to reproduce the proof written in the case $d=2$. Note that 
\begin{equation*}
	|G_{\omega^{(j,0)}}\triangle G_{\omega^{(j,1)}}| = \left|(G_{\omega^{(j,0)}}\cap A_j)\triangle (G_{\omega^{(j,1)}}\cap A_j)\right|, 
\end{equation*}
for all $\omega\in\{0,1\}^{M_\eta}$ and $j\in\{1\ldots,M_\eta\}$, and this quantity is equal to 
\begin{equation*}
	\int_{0}^\frac{\eta^2}{4} |B_{d-1}(0,\sqrt{r-r^2})|_{d-1}dr,
\end{equation*}
since as mentioned before $\eta^2/4$ is the height of the cap cut off by $H_j$, or in order words the distance between $y_j$ and the hyperplane $H_j$, independent of the index $j$. Therefore, 
\begin{align*}
	|G_{\omega^{(j,0)}}\triangle G_{\omega^{(j,1)}}| & = \int_{0}^\frac{\eta^2}{4} |B_{d-1}(0,\sqrt{r-r^2})|_{d-1}dr \\
	& = \int_{0}^\frac{\eta^2}{4} \beta_{d-1}(r-r^2)^{(d-1)/2}dr \\
	& = \beta_{d-1}\int_0^\frac{\eta^2}{4} (r-r^2)^{(d-1)/2}dr \\
	& = \frac{\beta_{d-1}\eta^{d+1}}{4^{d+1}}\int_0^1 u^{(d-1)/2}\left(1-\frac{\eta^2 u}{4}\right)^{(d-1)/2}du.
\end{align*}

Since $0<\eta^2/4<1/4$, we then get
\begin{equation}
	\label{1432}\frac{3^{(d-1)/2}\eta^{d+1} \beta_{d-1}}{2^{3d}(d+1)} \leq |G_{\omega^{(j,0)}}\triangle G_{\omega^{(j,1)}}| \leq \frac{\eta^{d+1} \beta_{d-1}}{2^{2d+1}(d+1)}.
\end{equation}
Now, continuing \eqref{25121986} and \eqref{25061986}, replacing $M$ by $M_\eta$ and $h$ by the lower bound in \eqref{1432} and using lemmas 3 and 4, we get that
\begin{equation}
	\label{1434}\sup_{G\in\mathcal C_d}\mathbb E_G\left[|G\triangle \hat G|\right] \geq C_8\eta^2\left(1-C_9\eta^{d+1} \right)^{2n},
\end{equation}
where 
\begin{equation*}
	C_8=\frac{3^{(d-1)/2}\beta_{d-1}d}{2^{4d+1}(d+1)\sqrt{d+2}}
\end{equation*}
and 
\begin{equation*}
	C_9=\frac{(1-e^{-\frac{1}{8\sigma^2}})\beta_{d-1}}{2^{2d+1}(d+1)}.
\end{equation*}
Note that since the ball $B_{d-1}(0,1/2)$ is included in the $(d-1)$-dimensional hypercube centered at the origin, with sides of length $1$, the following inequality holds
\begin{equation*}
	|B_{d-1}(0,\frac{1}{2})|=\frac{\beta_{d-1}}{2^{d-1}}\leq 1,
\end{equation*}
and this shows that $C_9<1$. Therefore, since $\eta<1$ as well, the concavity of the logarithm leads \eqref{1434} to
\begin{equation*}
	\sup_{G\in\mathcal C_d}\mathbb E_G\left[|G\triangle \hat G|\right] \geq C_8\eta^2\exp\left(2n\ln(1-C_9)\eta^{d+1}\right).
\end{equation*}
Let us choose $\eta=n^{-1/(d+1)}$, so that \eqref{1434} becomes 
\begin{equation*}
	\sup_{G\in\mathcal C_d}\mathbb E_G\left[|G\triangle \hat G|\right] \geq C_{10}n^{-\frac{2}{d+1}},
\end{equation*}
where $C_{10}=C_8(1-C_9)^2>0$.
$\blacksquare$

\paragraph{Proof of Theorem 5}

Let $r^*$ be a given and finite integer such that $d+1\leq r^*\leq R_n$. Note that if $r^*\leq r\leq r'$, then $\mathcal P_{r^*}\subseteq\mathcal P_{r}\subseteq\mathcal P_{r'}$. Therefore if $P\in\mathcal P_{r^*}$ and $G=P$ in model \eqref{posprob}, by Theorem 1 it is likely that with high probability we have, using the triangle inequality :
\begin{equation}
	\label{p_adapt1}|\hat P_n^{(r)}\triangle \hat P_n^{(r')}|\leq \frac{Cdr'\ln n}{n}, 
\end{equation}
for any $r^*\leq r\leq r'$, where $C$ is a constant. Therefore it is reasonable to select $\hat r$ as the minimal integer that satisfies \eqref{p_adapt1}.
 
Let $\hat r$ be chosen as in Theorem 5. For $r=d+1,\ldots,R_n$, let us denote by $A_r$ the event :
\begin{equation*} 
A_r = \left\{\forall r'=r,\ldots,R_n, |\hat P_n^{(r)}\triangle\hat P_n^{(r')}|\leq \frac{6dr'\ln n}{C_2n}\right\},
\end{equation*}
where $C_2$ is the same constant as in Theorem 1. Then $\hat r$ is the smallest integer $r\leq R_n$ such that $A_r$ holds.

Let $P\in\mathcal P_{r^*}$. We write the following :
\begin{equation}
	\label{p_adapt2} \mathbb E_P[|\hat P_n^{adapt}\triangle P|]  =\mathbb E_P[|\hat P_n^{adapt}\triangle P|I(\hat r\leq r^*)] +  \mathbb E_P[|\hat P_n^{adapt}\triangle P|I(\hat r> r^*)],
\end{equation}
and we bound separately the two terms in the right side. Note that if $\hat r\leq r^*$, then, since the event $A_{\hat r}$ holds by definition, 
\begin{equation*}
	|\hat P_n^{(r^*)}\triangle\hat P_n^{(\hat r)}|\leq \frac{6dr^*\ln n}{C_2n}.
\end{equation*}
Therefore, using the triangle inequality,
\begin{align}
	\mathbb E_P[|\hat P_n^{adapt}\triangle P|I(\hat r\leq r^*)] & \leq \mathbb E_P[|\hat P_n^{adapt}\triangle \hat P_n^{(r^*)}|I(\hat r\leq r^*)] + \mathbb E_P[|\hat P_n^{(r^*)}\triangle P|I(\hat r\leq r^*)] \nonumber \\
	& \leq \frac{6dr^*\ln n}{C_2n}+\frac{A_1dr^*\ln n}{n} \mbox{ }\mbox{by Corollary 1} \nonumber \\
	\label{p_adapt5}& \leq \frac{C_{11}r^*\ln n}{n},
\end{align}
where $C_{11}$ depends only on $d$ and $\sigma$.
The second term of \eqref{p_adapt2} is bounded differently. First note that for all $r=d+1,\ldots,R_n$, $\hat P_n^{(r)}\subseteq[0,1]^d$, so $|\hat P_n^{(r)}|\leq 1$. Thus, if $\overline {A_{r^*}}$ stands for the complement of the event $A_{r^*}$, we have the following inequalities.

\begin{align}
	 \mathbb E_P[|\hat P_n^{adapt}\triangle P|I(\hat r> r^*)] & \leq 2\mathbb P_P[\hat r> r^*] \nonumber\\
	 & \leq 2\mathbb P_P\left[\overline{A_{r^*}}\right] \nonumber \\
	 & \leq 2\sum_{r=r^*}^{R_n}\mathbb P_P\left[|\hat P_n^{(r^*)}\triangle\hat P_n^{(r)}|> \frac{6dr\ln n}{C_2n}\right] \nonumber
\end{align}
\begin{align}
	 & \leq 2\sum_{r=r^*}^{R_n}\mathbb P_P\left[|\hat P_n^{(r^*)}\triangle P|+|\hat P_n^{(r)}\triangle P|> \frac{6dr\ln n}{C_2n}\right] \nonumber\\
	 & \leq 2\sum_{r=r^*}^{R_n}\left (\mathbb P_P\left[|\hat P_n^{(r^*)}\triangle P|> \frac{3dr\ln n}{C_2n}\right]+\mathbb P_P\left[|\hat P_n^{(r)}\triangle P|> \frac{3dr\ln n}{C_2n}\right]\right) \nonumber \\
	 \label{p_adapt3}& \leq 2\sum_{r=r^*}^{R_n}\left (\mathbb P_P\left[|\hat P_n^{(r^*)}\triangle P|> \frac{3dr^*\ln n}{C_2n}\right]+\mathbb P_P\left[|\hat P_n^{(r)}\triangle P|> \frac{3dr\ln n}{C_2n}\right]\right).
\end{align}
Note that since $P\in\mathcal P_{r^*}$, it is also true that $P\in\mathcal P_{r}, \forall r\geq r^*$. Therefore, by Theorem 1, using first $x=\frac{dr^*\ln n}{C_2}$, then $x=\frac{dr\ln n}{C_2}$, it comes from \eqref{p_adapt3} that :
\begin{align}
	\mathbb E_P[|\hat P_n^{adapt}\triangle P|I(\hat r> r^*)] & \leq 2\sum_{r=r^*}^{R_n}\left (  C_1e^{-dr^*\ln n} + C_1e^{-dr\ln n} \right) \nonumber \\
	& \leq 4C_1R_nn^{-dr^*} \nonumber \\
	& \leq 4C_1R_nn^{-d(d+1)} \nonumber \\
	\label{p_adapt4} & \leq 4C_1\left(\frac{n}{\ln n}\right)^{\frac{d-1}{d+1}}n^{d(d+1)}.
\end{align}
Finally, using \eqref{p_adapt5} and \eqref{p_adapt4},
\begin{equation*}
	\mathbb E_P[|\hat P_n^{adapt}\triangle P|] \leq \frac{C_{12}r^*\ln n}{n},
\end{equation*}
where $C_{12}$ is a positive constant that depends on $d$ and $\sigma$. 
Let us now assume that $r^*$ is a given integer larger than $R_n$, possibly infinite, and that $P\in\mathcal P_{r^*}$. As in Theorem 5, if $r^*=\infty$ we denote by $\mathcal P_\infty$ the class $\mathcal C_d$.
Then with probability one, $\hat r\leq r^*$. First of all, note that obviously, since by definition, $\hat r\leq R_n$,
\begin{equation*}
	|\hat P_n^{(R_n)}\triangle\hat P_n^{(\hat r)}|\leq \frac{6dR_n\ln n}{C_2n}\leq \frac{6d}{C_2}\left(\frac{\ln n}{n}\right)^{\frac{2}{d+1}}.
\end{equation*}
Then, by the triangle inequality,
\begin{align*}
	\mathbb E_P[|\hat P_n^{adapt}\triangle P|] & \leq \frac{6d}{C_2}\left(\frac{\ln n}{n}\right)^{\frac{2}{d+1}} + \mathbb E_P[|\hat P_n^{(R_n)}\triangle P|] \\
	& \leq  \frac{6d}{C_2}\left(\frac{\ln n}{n}\right)^{\frac{2}{d+1}} + A'_1\left(\frac{\ln n}{n}\right)^\frac{2}{d+1},
\end{align*}
by Corollary 2, since $P\in\mathcal P_{r^*}\subseteq\mathcal P_{\infty}$ and $\hat P_n^{(R_n)}$ is the estimator of Theorem 3.
Theorem 5 is then proven. 
$\square$

\section{Appendix : proof of the lemmas}

\paragraph{Proof of Lemma 1}
Let us first state the following lemma, which gives the Steiner formula in the case of polytopes. It can also be found in \cite{3'}.
\begin{lemma}
For any polytope $R\subseteq \mathbb R^d$ the volume of $R^\lambda$ is polynomial in $\lambda$, with degree $d$, that is there exists $(L_0(R), \ldots, L_d(R)) \in \mathbb R^{d+1}$
\begin{equation*}
	|R^\lambda|=\sum_{k=0}^d L_k(R)\lambda^k, \mbox{   }\forall \lambda\geq0.
\end{equation*}
Besides, $L_0(R)=|R|$, $L_1(R)$ is the surface area of $R$ and $L_d(R)=|B_d(0,1)|$, independent of $R$, and all the $L_i(R), i=0,\ldots,d$ are nonnegative.
\end{lemma}

Note that in this lemma, if $R$ is included in $B_d(a,u)$ for some $a\in\mathbb R^d$ and $u>0$, then for all positive $\lambda$, 
\begin{equation*}
	R^\lambda\subseteq B_d(a,u)^\lambda=B_d(a,u+\lambda)
\end{equation*}
and if we denote by $\beta_d=|B_d(0,1)|$,
\begin{equation}
	\label{802479}|R^\lambda|=\sum_{k=0}^d L_k(R)\lambda^k \leq (u+\lambda)^d\beta_d.
\end{equation}
Therefore, since all the $L_i(R)$ are nonnegative, one gets  
\begin{equation}
	\label{SteinerUB} L_i(R)\leq(u+1)^d\beta_d, i=1,\ldots, d
\end{equation}
by taking $\lambda=1$ in \eqref{802479}.

Let $r\leq n$, and $P\in \mathcal P_r$. The polytope $P^*$ is constructed as follows. For any vertex $x$ of $P$, let $x^*$ be the closest point to $x$ in $[0,1]^d$ with coordinates that are integer multiples of $\frac{1}{n}$ (if there are several such points $x^*$, then one can take any of them). The euclidean distance between $x$ and $x^*$ is bounded by $\frac{\sqrt d}{n}$. 

Let us define $P^*$ as the convex hull of all these resulting $x^*$. Then $P^*\in\mathcal P_r^{(n)}$. 

For any set $G\subseteq\mathbb R^d$ and $\epsilon>0$ we denote by $G^\epsilon$ the set  
\begin{equation*}
	G^\epsilon=G+\epsilon B_d(0,1)=\{x\in\mathbb R^d : \rho(x,G)\leq \epsilon\}.
\end{equation*}
It is clear that the Hausdorff distance between $P$ and $P^*$ is less than $\frac{\sqrt d}{n}$. Therefore if we denote $\epsilon=\frac{\sqrt d}{n}$ we have $P^*\subseteq P^\epsilon$ and $P\subseteq (P^*)^{\epsilon}$. 

Since the two polytopes $P$ and $P^*$ are included in $B_d\left(a,\frac{\sqrt d}{2}\right)$, for $a=(1/2,\ldots,1/2)$, one gets from \eqref{SteinerUB} that 
\begin{equation*}
	L_i(R)\leq\left(\frac{\sqrt d}{2}+1\right)^d\beta_d \leq \left(\frac{3\sqrt d}{2}\right)^d, \mbox{    } i=0,\ldots,d
\end{equation*}
for $R=P$ or $P^*$.

We can now bound the Nikodym distance between $P$ and $P^*$  
\begin{align*}
	|P\triangle P^*| = |P\backslash P^*|+|P^*\backslash P| & \leq |(P^*)^\epsilon \backslash P^*| + |P^\epsilon\backslash P| \\
	& \leq 2\left(\frac{3\sqrt d}{2}\right)^d\beta_d\sum_{k=1}^d\left(\frac{\sqrt d}{n}\right)^k \leq \frac{2d^{d+1}(3/2)^d\beta_d}{n}.
\end{align*}
$\square$

\paragraph{Proof of Lemma 3}
First note that if $G\subseteq[0,1]^d$, then the density of the probability measure $\mathbb P_G$ with respect to the Lebesgue measure on $[0,1]^d\times \mathbb R$ is 
\begin{equation*}
	p_G(x,y)=\frac{1}{\sqrt{2\pi\sigma^2}}e^{-\frac{1}{2\sigma^2}(y-I(x\in G))^2}.
\end{equation*}
Therefore, by a simple algebra, if $G_1$ and $G_2$ are two subsets of $[0,1]^d$, then
\begin{align*}
	& \int_{[0,1]^2\times\mathbb R}\sqrt{p_{G_1}(x,y)p_{G_2}(x,y)}dxdy \\
	& = \int_{[0,1]^2}\exp\left(-\frac{I(x\in G_1\triangle G_2)}{8\sigma^2}\right)dx \\
	& = |G_1\triangle G_2|e^{-\frac{1}{8\sigma^2}}+1-|G_1\triangle G_2|, 
\end{align*}
and Lemma 3 follows from \cite{30}, Section 2.4.
$\square$

\paragraph{Proof of Lemma 4}

The fact that any $\eta$-packing family of $(S,\rho)$ is finite is clear and comes from the fact that $S$ is compact. Consider now a maximal $\eta$-packing family of $(S,\rho)$, denoted by $\{y_1,\ldots,y_{M_\eta}\}$. The surface area of $B_d(y_j,\eta/2)\cap S$ is independent of $j\in\{1,\ldots,M_\eta\}$, and we denote it by $V(\eta/2)$. A simple application of the Pythagorean theorem shows that $B_d(y_j,\eta/2)\cap S$ is a cap of height $\eta^2/4$ of $S$. Therefore, using Lemma 2.3 of \cite{25'}
\begin{equation*}
	V(\eta/2)\geq \beta_{d-1}\left(1-\frac{\eta^2}{4}\right)^{(d-3)/2}\eta^{d-1}.
\end{equation*}
Besides, since $\{y_1,\ldots,y_{M_\eta}\}$ is an $\eta$-packing family , the caps $B_d(y_j,\eta/2)\cap S, j=1,\ldots,M_\eta$ are pairwise disjoint and the surface area of their union is less than the surface area of $S$, which is equal to $\displaystyle{\frac{d\beta_d}{2^{d-1}}}$, so we get 
\begin{equation*}
	M_\eta V(\eta/2) \leq \frac{d\beta_d}{2^{d-1}}.
\end{equation*}
Therefore,
\begin{equation*}
	M_\eta \leq \frac{d\beta_d}{2^{d-1}V(\eta/2)} \leq \frac{d\beta_d}{2^{d-1}\beta_{d-1}\left(1-\frac{\eta^2}{4}\right)^{(d-3)/2}\eta^{d-1}}.
\end{equation*}
and the right inequality of Lemma 4 follows from the fact that $\eta^2/4\leq1/4$ and Lemma 2.2 of \cite{25'} which states that
\begin{equation}
	\label{Lemma2.2of25'} \frac{\sqrt{2\pi}}{\sqrt{d+2}} \leq \frac{\beta_d}{\beta_{d-1}} \leq \frac{\sqrt{2\pi}}{\sqrt{d}}.
\end{equation}
The left inequality of Lemma 4 comes from the fact that any maximal $\eta$-packing family is an $\eta$-net. Indeed, consider a maximal $\eta$-packing family $\mathcal Y$, and assume it is not an $\eta$-net. Then there exists $x\in S$ such that for all $y\in\mathcal Y$, $\rho(x,y)>\epsilon$. Therefore $\{x\}\cup\mathcal Y$ is an $\eta$-net that contains $\mathcal Y$ strictly. This contradicts maximality of $\mathcal Y$. Therefore the family $\{y_1,\ldots,y_{M_\eta}\}$ is an $\eta$-net of $S$, and the caps $B_d(y_j,\eta)\cap S, j=1,\ldots,M_\eta$ cover the sphere $S$, so that
\begin{equation*}
	M_\eta V(\eta) \geq \frac{d\beta_d}{2^{d-1}}.
\end{equation*}
Using again Lemma 2.3 of \cite{25'}, we bound $V(\eta)$ from above 
\begin{equation*}
	V(\eta)\leq \beta_{d-1}\eta^{d-1},
\end{equation*}
and then the desired result follows again from \eqref{Lemma2.2of25'}. 
$\square$

\selectlanguage{english}

\end{document}